\input amssym
\font\tiny=cmr8
\def\Q{{\Bbb{Q}}}
\def\lie{{\frak{lie}}}
\def\krv{{\frak{krv}}}
\def\der{{\frak{der}}}
\def\grt{{\frak{grt}}}
\def\ds{{\frak{ds}}}

\def\krv{{\frak{krv}}}
\centerline{\bf On the elliptic Kashiwara-Vergne Lie algebra}
\vskip .5cm
\centerline{\'Elise Raphael and Leila Schneps\footnote{$^*$}{Both authors
gratefully acknowledge support from the Simons Center for Geometry and
Physics, Stony Brook University, at which some of the research for this paper
was performed.}}
\vskip .5cm
{\narrower{\narrower{\noindent{\bf Abstract.} We recall the definitions of two 
independently defined elliptic versions of the Kashiwara-Vergne Lie algebra 
$\krv$, namely the Lie algebra $\krv^{(1,1)}$ constructed by A.~Alekseev, 
N.~Kawazumi, Y.~Kuno and F.~Naef arising from the study of graded formality 
isomorphisms associated to topological fundamental groups of surfaces,
and the Lie algebra $\krv_{ell}$ defined using mould theoretic techniques
arising from multiple zeta theory by E.~Raphael and L.~Schneps, and show that 
they coincide.}\par}\par}
\vskip .5cm
\noindent {\bf \S 1. Introduction}
\vskip .2cm
From its inception in Grothendieck's {\it Esquisse d'un Programme} ([G]), 
Grothendieck-Teichm\"uller theory was intended to study the automorphism
groups of the profinite mapping class groups -- the
fundamental groups of moduli spaces of Riemann surfaces
of all genera and any number of marked points -- with the goal of 
discovering new properties of the absolute Galois group ${\rm Gal}(\overline{\Bbb Q}/{\Bbb Q})$.  However, due to the ease
of study of the genus zero mapping class groups, which are essentially
braid groups, the genus zero case garnered most of the attention, 
starting from the definition of the Grothendieck-Teichm\"uller group 
$\widehat{GT}$ by V.G.~Drinfel'd ([D]) and the simultaneous construction by
Y.~Ihara of the Grothendieck-Teichm\"uller Lie algebra $\grt$ ([I1,I2]) 
in 1991. The extension of the definition to a Grothendieck-Teichm\"uller group 
acting on the profinite mapping class groups in all genera was subsequently 
given in 2000 by A.~Hatcher, P.~Lochak, L.~Schneps and H.~Nakamura (cf.~[HLS], 
[NS]).  The higher genus profinite Grothendieck-Teichm\"uller group satisfies 
the two-level principle articulated by Grothendieck, which states that
the subgroup of $\widehat{GT}$ consisting of automorphisms that extend to
the genus one mapping class groups with one or two marked points will 
automatically extend to automorphisms of the higher mapping class groups.

It has proved much more difficult to extend the Lie algebra Grothendieck-Teichm\"uller construction to higher genus.
Indeed, while the genus zero mapping class groups have a natural Lie algebra 
analog in the form of the braid Lie algebras, there is
no good Lie algebra analog of the higher genus mapping class groups.  The
only possible approach for the moment seems to be to replace the higher
genus mapping class groups by their higher genus braid subgroups, which do 
have good Lie algebra analogs\footnote{$^{**}$}{Another approach would be to
replace the higher genus mapping class groups by their Torelli subgroups, 
which also have good Lie algebraic analogs determined by R.~Hain [H]. 
In particular, this would include the key case of higher genus with $0$ marked 
points, which have no associated braid groups.  However, there has been no 
development of Lie Grothendieck-Teichm\"uller theory in this context as yet.}.
An early piece of work due to H.~Tsunogai ([T]) in 2003
computed the relations that must be satisfied by a derivation acting on
the genus one 1-strand braid Lie algebra (which is free on two generators)
to ensure that it extends to a derivation on the genus one 2-strand braid
Lie algebra, in analogy with the derivations in $\grt$, shown by Ihara to
be exactly those that act on the genus zero 4-strand braid Lie algebra (also a 
free Lie algebra on two generators) and extend to derivations of the 
5-strand braid Lie algebra.

After this, the next real breakthrough in the higher genus Lie algebra 
situation came with the work of B.~Enriquez ([En], 2014) based on his
previous joint work with D.~Calaque, P.~Etingof ([CEE], 2009). 
In particular, using the same approach as 
Tsunogai of replacing the higher genus mapping class groups with their higher 
genus braid subgroups, Enriquez in [En] 
was able to extend the definition of $\grt$ to 
an elliptic 
version $\grt_{ell}$, which he identified with an explicit Lie subalgebra
of the algebra of derivations of the algebra of the genus one 1-strand 
braid Lie algebra that extend to derivations of the 2-strand genus one
braid Lie algebra.  He showed in particular
that there is a canonical surjection $\grt_{ell}\rightarrow\!\!\!\!\rightarrow
\grt$, and a canonical section of this surjection, $\gamma:\grt\hookrightarrow
\grt_{ell}$.
\vskip .2cm
The Grothendieck-Teichm\"uller Lie algebra is closely related
to two other Lie algebras, the double shuffle Lie algebra that arises from
the theory of multiple zeta values and the Kashiwara-Vergne Lie algebra that
arises from solutions to the (linearized) Kashiwara-Vergne problem.  Indeed,
there exist injective Lie algebra morphisms
$$\grt\hookrightarrow\ds\hookrightarrow\krv,$$
by work of H.~Furusho ([F]) for the first injection, J.~\'Ecalle and
L.~Schneps ([E],[S1]) for the second and A.~Alekseev and C.~Torossian ([AT])
for a direct proof that $\grt$ maps into $\krv$.  In fact, these three algebras
are conjectured to be isomorphic, a conjecture that is supported by 
computation of the graded parts up to weight about 20.  Thus it was a
natural consequence of the work of Enriquez to consider the possibility
of extending also these other Lie algebras from genus zero to genus one.  
An answer was proposed for the double shuffle Lie algebra in [S3],
which proposes a definition of an elliptic double shuffle Lie algebra
$\ds_{ell}$ based on mould theory and an elliptic interpretation of a strong 
theorem due to J.~\'Ecalle (cf.~[E],[S2],[S3]).  The elliptic double shuffle
Lie algebra is compatible with Enriquez's construction; in particular
there is an injective Lie algebra morphism $\gamma:\ds\hookrightarrow \ds_{ell}$
which extends Enriquez's section in the sense that the image of the
Lie subalgebra $\grt\subset \ds$ is equal to Enriquez's image $\gamma(\grt)$.

One interesting aspect of the mould theoretic approach is that it reveals a 
close relationship between the elliptic double shuffle Lie algebra and the 
associated graded of the usual double shuffle Lie algebra for the depth 
filtration.  In the article [RS], the authors of this paper showed that an 
analogous approach works to construct an elliptic version of $\krv$, denoted 
$\krv_{ell}$, which is given by two defining mould theoretic properties, 
and again has the key features of 
\vskip .2cm
\noindent $\bullet$ being naturally identified with a Lie subalgebra of
the derivation algebra of the free Lie algebra on two generators; 
\vskip .1cm
\noindent $\bullet$ being equipped with an injective Lie algebra morphism 
$\gamma:\krv\hookrightarrow\krv_{ell}$ which extends the 
Grothendieck-Teichm\"uller and double shuffle maps;
\vskip .1cm
\noindent $\bullet$ having a structure closely related to that of 
the associated graded of $\krv$ for the depth filtration.
\vskip .2cm
In independent work, A.~Alekseev et al.~([AKKN1,2]) took a different approach 
to the construction of higher genus Kashiwara-Vergne Lie algebras 
$\krv^{(g,n)}$ for all $g,n\ge 1$, following the classical approach to the
Kashiwara-Vergne problem which focuses on determining graded formality
isomorphisms between prounipotent fundamental groups of surfaces and
their graded counterparts (i.e.~the exponentials of the associated gradeds
of their associated Lie algebras).

More precisely, if $\Sigma$ denotes a compact oriented surface of genus $g$ 
with $n + 1$ 
boundary components, the space $g(\Sigma)$ spanned by free homotopy classes of 
loops in $\Sigma$ carries the structure of a Lie bialgebra equipped with the 
Goldman bracket and the Turaev cobracket. The Goldman-Turaev formality 
problem is the construction a Lie bialgebra homomorphism $\theta$ from 
$g(\Sigma)$ to its associated graded ${\rm gr}\,g(\Sigma)$ such that 
${\rm gr}\,\theta = {\rm id}$.  In order to solve this problem, Alekseev et
al.~defined 
a family $KV(g, n + 1)$ of Kashiwara-Vergne problems.  In the particular
situation where $(g,n)=(1,0)$, the surface $\Sigma$ is of genus 1 with one boundary 
component, and its fundamental group is free on two generators $A$, $B$,
with the boundary loop being given by $C=(A,B)$.  The prounipotent fundamental
group is then free on two generators $e^a$, $e^b$ with a boundary element
$e^c$ satisfying $e^c=(e^a,e^b)$.  The associated Lie algebra is free on
generators $a$, $b$ and the logarithm of the boundary loop is given by
$$c=CH\Bigl(CH\bigl(CH(a,b),-a\bigr),-b\Bigr)=[a,b]+{\rm higher\ order\ terms},
$$
where $CH$ denotes the Campbell-Hausdorff law on ${\rm Lie}[a,b]$.  To
define the genus one Kashiwara-Vergne Lie algebra $\krv^{(1,1)}$, Alekseev
et al.~first defined the space of
derivations $u$ of ${\rm Lie}[a,b]$ that annihilate the boundary element $c$
and further satisfy a certain non-commutative divergence condition
(see \S 2 for more detail), and then took $\krv^{(1,1)}$ to be the associated
graded of the above space, which comes down to using the same defining
conditions with $c$ replaced by $[a,b]$.  They showed that the resulting
space is a Lie algebra under the bracket of derivations, and also 
that, like $\krv_{ell}$, it is equipped with an injective Lie algebra morphism 
$\krv\hookrightarrow \krv^{(1,1)}$ that extends the morphism
$\grt\hookrightarrow \grt_{ell}$ constructed by Enriquez.
\vskip .2cm
The main result of this article is the equivalence of these two definitions
of the elliptic Kashiwara-Vergne Lie algebra.
\vskip .2cm
\noindent {\bf Main Theorem.}  {\it There is a canonical isomorphism
$\krv^{(1,1)}\simeq \krv_{ell}.$}
\vskip .2cm
It is an easy consequence of known results that the first defining property
of $\krv_{ell}$ corresponds to the annihilation of $c$.  The proof of the 
theorem thus consists essentially in comparing the second defining properties
of the two algebras.  The article is organised as follows.
In \S 2, we recall the definition of $\krv^{(1,1)}$, 
in particular the divergence property, and in \S 3, we give a new
reformulation of the divergence property. In \S 4 we recall the definition
of $\krv_{ell}$ and show that translating its second 
mould theoretic defining property back to a property of
derivations on ${\rm Lie}[a,b]$, it coincides with the reformulated version
of the divergence property of $\krv^{(1,1)}$ given in \S 3, which completes
the proof.
\vskip .2cm
\noindent {\bf \S 2. The elliptic Kashiwara-Vergne Lie algebra from [AKKN].}
Let $\lie^{(1,1)}$ be the free Lie algebra on two generators ${\rm Lie}[a,b]$,
to be thought of as the Lie algebra associated to the fundamental group of
the once-punctured torus.  We set $c=[a,b]$ so that the relation
$[a,b]=c$  holds in $\lie^{(1,1)}$. Let $\lie_n^{(1,1)}$ denote the
weight $n$ part of $\lie^{(1,1)}$, where the weight is the total degree in
$a$ and $b$, and let $\lie_{n,r}^{(1,1)}$ denote the weight $n$, 
depth $r$ part of $\lie^{(1,1)}$, where the depth is the $b$-degree.

For any element $f\in \lie^{(1,1)}$, we decompose it as
$$f=f_aa+f_bb=af^a+bf^b=af^a_aa+af^a_bb+bf^b_aa+bf^b_bb.$$
Let $\der^{(1,1)}$ denote the Lie subalgebra 
of ${\rm Der}\,\lie^{(1,1)}$ of derivations $u$ such that $u(c)=0$. Let
$\der_n^{(1,1)}$ denote the subspace of $\der^{(1,1)}$ of derivations $u$
such that $u(a),u(b)\in \lie_n^{(1,1)}$.

Let $f\in \lie_n^{(1,1)}$ for $n>1$.
Then $f$ is the value on $a$ of a derivation $u\in \der_n^{(1,1)}$ if and only 
if $f$ is push-invariant, in which case $u(b)$ is uniquely defined.  If 
$u\in\der^{(1,1)}$ and $u(a)=f$, $u(b)=g$, we sometimes write $u=D_{f,g}$.

Let $Tr_2$ be the quotient of $\lie^{(1,1)}$ by the equivalence relation:
two words $w$ and $w'$ are equivalent if one can be obtained from the other
by cyclic permutation of the letters.

The {\it elliptic divergence} $div:\der^{(1,1)}\rightarrow Tr_2$ is
defined by
$$div(u)=tr(f_a+g_b)$$
where $u=D_{f,g}$.  Since $u([a,b])=[a,g]+[f,b]=0$, we have
$$ag_aa+ag_bb-ag^aa-bg^bb=bf_aa+bf_bb-af^ab-bf^bb.$$
Comparing the terms on both sides that start with $a$ and end with $b$ shows
that $g_b=-f^a$.  Thus we can write the divergence condition as a function
of just $f$:
$$div(u)=tr(f_a-f^a).$$
\vskip .3cm
\noindent {\bf Definition.} The elliptic Kashiwara-Vergne Lie algebra 
${\krv}^{(1,1)}$ defined in [AKKN] is the $\Q$-vector space spanned by
the derivations $u\in\der_n^{(1,1)}$, $n\ge 3$, such that
$$div(u)= \cases{K\,tr\bigl(c^{{n-1}\over{2}}\bigr)\ \hbox{for some}\ 
K\in \Q&if $n$ is odd\cr
0&if $n$ is even.}\eqno(1)$$
It is closed under the bracket of derivations.
\vskip .3cm
\noindent {\bf \S 3. A reformulation of the $div$ condition.} Observe that
the Lie algebra $\krv^{(1,1)}$ is bigraded by the weight and depth. We
write $\krv_{n,r}^{(1,1)}$ for the vector subspace of $\krv_{ell}$
spanned by derivations $u$ such that $u(a)\in\lie_{n,r}^{(1,1)}$ and
$u(b)\in\lie_{n,r+1}^{(1,1)}$. Let $u\in \krv_{n,r}^{(1,1)}$, let
$f=u(a)$ and $g=u(b)$, and write
$$f = \sum_{\underline{i}=(i_0,...,i_r)} c_{\underline{i}}\ a^{i_0}ba^{i_1}b \cdots ba^{i_r}.$$
Then 
$$f_a =\sum_{\underline{i}\ {\rm s.t.}\ i_r\ne 0} 
c_{\underline{i}}\ a^{i_0}ba^{i_1}b \cdots ba^{i_r-1}$$ 
and
$$f^a =\sum_{\underline{i}\ {\rm s.t.}\ i_0\ne 0} 
c_{\underline{i}}\ a^{i_0-1}ba^{i_1}b \cdots ba^{i_r}.$$ 

For any word $w$ in $a,b$, let $C(w)$ denote the trace class of $w$,
i.e. the set of words obtained by cyclically permuting the letters of $w$.
The trace of a polynomial $h$ is given by
$$\sum_{C(w)} \bigl(tr(h)\,|\,C(w)\bigr)\cdot C(w)$$
where the sum runs over the trace classes of words and the coefficient
$\bigl(tr(h)\,|\,C(w)\bigr)$ of the class $C(w)$ in $tr(h)$ is given by
$$\bigl(tr(h)\,|\,C(w)\bigr)=\sum_{v\in C(w)} (h|v),$$
where $(h|v)$ denotes the coefficient of $v$ in $h$.

Fix a word $w$ of weight $n-1$ and depth $r$ and consider the coefficient
of $C(w)$ in $tr(f_a-f^a)$:
$$\eqalign{\bigl(tr(f_a-f^a)\,|\,C(w)\bigr)
&=\sum_{v\in C(w)} (f_a|v)-(f^a|v)\cr
&=\sum_{v\in C(w)} (f_aa|va)-(af^a|av)\cr
&=\sum_{v\in C(w)} (f|va)-(f|av).}\eqno(2)$$
For all words $v\in C(w)$ that start in $a$, the word $v'\in C(w)$ obtained
from $v$ by taking the first $a$ of $v$ and putting it at the end satisfies
$va=av'$ and thus $(f|va)=(f|av')$.  Thus the corresponding terms in (2)
cancel out, i.e.~writing $C^a(w)$ (resp.~$C^b(w)$, $C_a(w)$, $C_b(w)$) for
the terms in $C(w)$ that start with $a$ (resp.~start with $b$, end with $a$,
end with $b$), we have
$$\sum_{v\in C^a(w)} (f|va)- \sum_{v\in C_a(w)} (f|av)=0,$$
so (2) reduces to
$$\bigl(tr(f_a-f^a)\,|\,C(w)\bigr)=\sum_{v\in C^b(w)} (f|va)- 
\sum_{v\in C_b(w)} (f|av).\eqno(3)$$
Since for any word $v$ of length $n-2$, if $v=bu\in C^b(w)$ then $ub\in C_b(w)$,
we can write (3) as
$$\bigl(tr(f_a-f^a)\,|\,C(w)\bigr)=\sum_{u\ {\rm s.t.}\ bu\in C^b(w)} (f|bua)- 
\sum_{u\ {\rm s.t.}\ ub\in C_b(w)} (f|aub).\eqno(4)$$

For any word $u=a^{i_0}b\cdots ba^{i_r}$ of depth $r$, we define 
$$push(u)=a^{i_r}ba^{i_0}b\cdots ba^{i_{r-1}}$$
and
$$pushsym(u)=\sum_{i=0}^{r} push^i(u),$$
and extend the operators $push$ and $pushsym$ to polynomials by linearity.
If $u$ is a word of weight $n-2$ and depth $r-1$ 
such that $bu\in C^b(w)$, then $ub\in C_b(w)$ and we have
$$\cases{C^b(w)=\{b\,push^i(u)\,|\,0\le i\le r-1\}\cr
C_b(w)=\{push^i(u)\,b\,|\,0\le i\le r-1\}.}$$
Note that $C^b(w)$ may be of order less than $r$ (in fact strictly
dividing $r$) when $w$ has a symmetry under the push.  In this case,
summing over the set of pushes of $u$ from $0$ to $r-1$ comes down
to summing $r/|C^b(w)|$ times over the set $C^b(w)$ or $C_b(w)$.
Using this, we rewrite (4) for $w=ub$ as
$$\eqalign{\sum_{bu\in C^b(w)} (f|bua)-\sum_{ub\in C_b(w)} (f|aub)&=
\sum_{bu\in C^b(w)} (f^b_a|u)-\sum_{ub\in C_b(w)} (f^a_b|u)\cr
&={{|C^b(w)|}\over{r}}\sum_{i=0}^{r-1} \bigl(f^b_a\,|\,push^i(u)\bigr)-
{{|C_b(w)|}\over{r}}\sum_{i=0}^{r-1} \bigl(f^a_b
\,|\, push^i(u)\bigr)\cr
&={{|C_b(w)|}\over{r}}\bigl(pushsym(f^b_a-f^a_b)\,|\,u\bigr).}$$
This allows us to rewrite the divergence condition (1) on an element
$f\in \krv_{n,r}^{(1,1)}$
as the following 
family of relations for all words $u$ of weight $n-2$ and depth $r-1$:
$$\bigl(pushsym(f^b_a-f^a_b)|u\bigr)=
\cases{{{Kr}\over{|C_b(ub)|}}\sum_{v\in C(ub)} \bigl([a,b]^r \,|\, v\bigr)
\ \hbox{for some}\ K\in \Q&if $n=2r+1$\cr
0&if $n\ne 2r+1$.}\eqno(5)$$
This is the version of the divergence condition that we will use for
comparison with the Lie algebra $\krv_{ell}$.
\vskip .3cm
\noindent{\bf \S 4. The mould theoretic $\krv_{ell}$ from [RS].}
Recall that a mould is a family $A=(A_r)_{r\ge 0}$ where $A_r(u_1,\ldots,u_r)$
is a function of $r$ commutative variables.  We restrict our attention here
to rational-function moulds with coefficients in $\Q$. These form a 
$\Q$-vector space under componentwise addition and multiplication by scalars. 
When the number of variables is specified, we drop the subscript $r$,
for instance we write $A(u_1,\ldots,u_r)=A_r(u_1,\ldots,u_r)$.

A mould is said to be {\it alternal} if $A(\emptyset)=0$ and
$$\sum_{w\in sh\bigl((u_1,\ldots,u_k),(u_{k+1},\ldots,u_r)\bigr)} A_r(w)=0$$
for $r\ge 2$ and $1\le k\le r-1$.  

Let us define a few mould operators.
The $swap$, $push$, $circ$ and $\Delta$-operators on moulds are given by
$$\eqalign{swap(A)(v_1,\ldots,v_r)&=A(v_r,v_{r-1}-v_r,\ldots,v_1-v_2)\cr
push(A)(u_1,\ldots,u_r)&=A(u_2,\ldots,u_r,-u_1-\cdots-u_r)\cr
circ(A)(v_1,\ldots,v_r)&=A(v_r,v_1,\ldots,v_{r-1})\cr
\Delta(A)(u_1,\ldots,u_r)&=(u_1+\cdots+u_r)u_1\cdots u_r\,
A(u_1,\ldots,u_r).}$$
There is no difference between the use of the commutative variables $u_i$
or $v_i$, however the $v_i$'s are traditionally used for operators and
relations concerning the swap of a mould. 

There is a direct connection between power series in $a,b$ (having no
constant term in $a$) and polynomial-valued moulds.  Let $c_i=ad(a)^{i-1}(b)$
for $i\ge 1$, and consider Lie algebra $Lie[c_1,c_2,\ldots]$
inside the polynomial algebra $\Q\langle c_1,c_2,\ldots\rangle$.
By Lazard elimination, these algebras are free and all polynomials in 
$\lie^{(1,1)}$ having no linear term in $a$ can be written as
Lie polynomials in the $c_i$.

There is a bijection between the space of polynomials
in the $c_i$ and the space of polynomial-valued moulds, coming
from linearly extending the map 
$$c_{a_1}\cdots c_{a_r}\rightarrow (-1)^{n+r}u_1^{a_1-1}\cdots u_r^{a_r-1},\eqno(6)$$
where $n=a_1+\cdots+a_r$. It is well-known that under this map,
the subspace ${\rm Lie}[c_1,c_2,\ldots]$ of $\lie^{(1,1)}$, which consists
of all Lie polynomials having no linear term in $a$, maps bijectively
onto the space of alternal polynomial-valued moulds.
Furthermore, if we write such an $f\in\lie^{(1,1)}$  as
$$f=\sum_{\underline{i}} c_{\underline{i}}\,a^{k_0}b\cdots ba^{k_r}$$
and $F$ for the corresponding mould, then $swap(F)$ is explicitly given by
$$swap(F)(v_1,\ldots,v_r)=\sum_{\underline{i}\ {\rm s.t.}\ k_r=0} c_{\underline{i}}\,v_1^{k_0}\cdots v_r^{k_{r-1}}.\eqno(7)$$

A mould $A$ is said to be {\it push-invariant} if $push(B)=B$, and
{\it circ-constant} if there exists $K\in \Q$ such that for all $r\ge 2$, we 
have
$$\sum_{i=0}^{r-1} circ^i(A)(v_1,\ldots,v_r)=Kr.$$
If $K=0$ then $A$ is said to be {\it circ-neutral}.
\vskip .2cm
\noindent{\bf Definition.} The mould version of $\krv_{ell}$ consists of
of all polynomial-valued moulds $F$ that are alternal and push-invariant
and such that $swap\bigl(\Delta^{-1}(F)\bigr)$ is circ-constant.  
\vskip .2cm
The space $\krv_{ell}$ is bigraded for the depth and the degree. Let
$F\in \krv_{ell}$ be a mould of depth $r$ and degree $d$, so that it
corresponds under the bijection (6) to a polynomial $f\in \lie_{n,r}^{(1,1)}$ with $n=d+r$.
It is well-known that the mould push-invariance property of a polynomial-valued
mould $F$ is equivalent to the polynomial push-invariance $push(f)=f$. 
In turn, the polynomial push-invariance of $f$ implies that there exists
a unique polynomial $g\in \lie_{n,r+1}^{(1,1)}$ such that
setting $u(a)=f$, $u(b)=g$, we obtain a derivation $u\in \der_n^{(1,1)}$. 
The Lie bracket on $\krv_{ell}$ corresponds to the Lie bracket on
$\krv^{(1,1)}$, namely bracketing of the derivations $u$.
Thus, in order to prove that $\krv_{ell}$ is in bijection with $\krv^{(1,1)}$, it remains
only to prove that the circ-constance condition on 
$\Delta^{-1}\bigl(swap(F)\bigr)$ is equivalent to the divergence condition
(5) on $f$.
\vskip .2cm
Since
$$\Delta^{-1}(F)(u_1,\ldots,u_r)={{1}\over{(u_1+\cdots+u_r)u_1\cdots u_r)}}
F(u_1,\ldots,u_r),$$
we have 
$$swap\bigl(\Delta^{-1}(F)\bigr)(v_1,\ldots,v_r)={{1}\over{v_1(v_1-v_2)\cdots
(v_{r-1}-v_r)v_r}}swap(F)(v_1,\ldots,v_r),$$
so the circ-constance condition is given explicitly by
$$ {{swap(F)(v_1,\dots,v_r)}\over{v_1(v_1-v_2)\dots(v_{r-1}-v_r)v_r}} + {{swap(F)(v_2,\dots,v_1)}\over{v_2(v_2-v_3)\dots(v_r-v_1)v_1}}+\dots + {{swap(F)(v_r,\dots,v_{r-1})}\over{v_r(v_r-v_1)\dots(v_{r-2}-v_{r-1})v_{r-1}}} = Kr$$
for all depths $r\ge 2$. 
Putting this over a common denominator gives the equivalent equality
$$swap(F)(v_1,v_2,\dots,v_r) v_2\dots v_{r-1} (v_r-v_1) +swap(F)(v_2,\ldots,v_r,v_1)v_3\cdots v_r(v_1-v_2)+\cdots$$
$$\ \ \ \ \ \ \ \ \ \ \ + swap(F)(v_r,\dots,v_{r-1}) v_1\dots v_{r-2} (v_{r-1}-v_r) 
= Kr v_1\cdots v_r(v_1-v_2)\cdots(v_r -v_1).\eqno(8)$$
The left-hand side of (8) expands to 
$$\eqalign{&v_2\dots v_{r-1}v_r
 \,swap(F)(v_1,\dots,v_r)-v_1v_2\dots v_{r-1}\,swap(F)(v_1,\dots,v_r)  \cr
&\ \ \ +v_1v_3\dots v_r\;swap(F)(v_2,\dots,v_r,v_1) -v_2v_3..v_r \,swap(F)(v_2,..,v_r,v_1)+\cdots \cr
&\ \ \ \ \ \ +v_1\dots v_{r-1}\; swap(F)(v_r,\dots,v_{r-1}) - v_1\dots v_{r-2}v_r\; swap(F)(v_r,\ldots v_{r-1}).}\eqno(9)$$
Fix a monomial $v_1^{i_1+1}v_2^{i_2+1}\dots v_r^{i_r+1}$. 
Calculating its coefficient in (9) yields
$$\eqalign{& \bigl(swap(F)(v_1,\dots,v_r) | v_1^{i_1+1}v_2^{i_2}\dots v_r^{i_r}\bigr) - \bigl(swap(F)(v_1,\dots,v_r) | v_1^{i_1}v_2^{i_2}\dots v_r^{i_r+1}\bigr)  \cr
&\ \ \ +\bigl(swap(F)(v_2,\dots,v_r,v_1) | v_1^{i_1}v_2^{i_2+1}\dots v_r^{i_r}\bigr) - \bigl(swap(F)(v_2,\dots,v_r,v_1) | v_1^{i_1+1}v_2^{i_2}\dots v_r^{i_r}\bigr)+\cdots\cr
&\ \ \ \ \ \ \ +\bigl(swap(F)(v_r,v_1,\dots,v_{r-1}) | v_1^{i_1}\dots v_r^{i_r+1}\bigr)-\bigl(swap(F)(v_r,v_1,\dots,v_{r-1}) | v_1^{i_1}\dots v_{r-1}^{i_{r-1}+1} v_r^{i_r}\bigr)\cr
&=\bigl(swap(F)(v_1,\dots,v_r) | v_1^{i_1+1}v_2^{i_2}\dots v_r^{i_r}\bigr) - \bigl(swap(F)(v_1,\dots,v_r) | v_1^{i_1}v_2^{i_2}\dots v_r^{i_r+1}\bigr)  \cr
&\ \ \ +\bigl(swap(F)(v_1,\dots,v_r) | v_1^{i_2+1}v_2^{i_3}\dots v_r^{i_1}\bigr) - \bigl(swap(F)(v_1,\dots,v_r) | v_1^{i_2}v_2^{i_3}\dots v_r^{i_1+1}\bigr)+\cdots\cr
&\ \ \ \ \ \ +\bigl(swap(F)(v_1,\dots,v_r) | v_1^{i_r+1}\dots v_r^{i_{r-1}}\bigr)- \bigl(swap(F)(v_r,v_1,\dots,v_{r-1}) | v_1^{i_r}\dots v_r^{i_{r-1}+1}\bigr),}$$
where the equality is obtained by bringing every term back to a coefficient
of a word in $swap(F)(v_1,\ldots,v_r)$.

The circ-constance condition on $swap\bigl(\Delta^{-1}(F)\bigr)$ can thus
be expressed by the family of relations for every tuple $(i_1,\ldots,i_r)$:
$$\eqalign{& \bigl(swap(F)(v_1,\dots,v_r) | v_1^{i_1+1}v_2^{i_2}\dots v_r^{i_r}\bigr) - \bigl(swap(F)(v_1,\dots,v_r) | v_1^{i_1}v_2^{i_2}\dots v_r^{i_r+1}\bigr)  \cr
&\ \ \ +\bigl(swap(F)(v_1,\dots,v_r) | v_1^{i_2+1}v_2^{i_3}\dots v_r^{i_1}\bigr) - \bigl(swap(F)(v_1,\dots,v_r) | v_1^{i_2}v_2^{i_3}\dots v_r^{i_1+1}\bigr)+\cdots\cr
&\ \ \ \ \ \ +\bigl(swap(F)(v_1,\dots,v_r) | v_1^{i_r+1}\dots v_r^{i_{r-1}}\bigr)- \bigl(swap(F)(v_r,v_1,\dots,v_{r-1}) | v_1^{i_r}\dots v_r^{i_{r-1}+1}\bigr)\cr
&=Kr\,\bigl((v_1-v_2)\dots (v_r-v_1)\,|\,v_1^{i_1}\dots v_r^{i_r}\bigr).
}\eqno(10)$$
We can use formula (7) to translate this equality directly back into
terms of the polynomial $f(a,b)$ corresponding to the mould $F$ under the
bijection (6).  We can write the right-hand side of (10) as
$$Kr\,\bigl((v_1-v_2)\cdots (v_{r-1}-v_r)v_r\,|\,v_1^{i_1}\cdots v_r^{i_r}\bigr)-
Kr\,\bigl((v_1-v_2)\cdots (v_{r-1}-v_r)v_1\,|\,v_1^{i_1}\cdots v_r^{i_r}\bigr),
\eqno(11)$$
or equivalently, setting
$B(v_1,\ldots,v_r)=(v_1-v_2)\cdots(v_{r-1}-v_r)v_r$, as
$$Kr\,\bigl(B\,|\,v_1^{i_1}\cdots v_r^{i_r}\bigr)-(-1)^{r-1}Kr\,\bigl(B\,|\,v_1^{i_r}\cdots v_r^{i_1}\bigr),\eqno(12)$$
by numbering the $v_i$ in the second term of (11) in the opposite order
\vskip .2cm

We have $[a,b]=ad(a)(b)=c_2$, so $[a,b]^r=c_2^r$,
and by formula (6), the polynomial-valued mould corresponding to 
$[a,b]^r$ is thus given by
$$A(u_1,\ldots,u_r)=(-1)^ru_1\cdots u_r.$$
The swap of this mould is given by
$$swap(A)(v_1,\ldots,v_r)=-(v_1-v_2)\cdots(v_{r-1}-v_r)v_r.$$
Thus the mould $B$ of (12) satisfies $B=-swap(A)$, 
so by (7), the expression (12) reformulating
the right-hand side of (10) translates back to polynomials as to
$$-Kr\,\Bigl([a,b]^r\,|\,a^{i_1}b\cdots a^{i_r}b\Bigr)+(-1)^{r-1}
Kr\,\Bigl([a,b]^r\,|\,a^{i_r}b\cdots a^{i_1}b\Bigr).\eqno(13)$$
Using (7) to directly translate the left-hand side of (10) in terms
of the polynomial $f$, we thus obtain the following expression equivalent
to the circ-neutrality property (10):
$$\eqalign{&\bigl(f\,|\,a^{i_1+1}ba^{i_2}b\dots ba^{i_r}b\bigr)-
\bigl(f\,|\,a^{i_1}ba^{i_2}b\cdots ba^{i_r+1}b\bigr)  \cr
&\ \ \ +\bigl(f\,|\,a^{i_2+1}ba^{i_3}b\dots ba^{i_1}b\bigr)-
\bigl(f\,|\,a^{i_2}ba^{i_3}b\cdots ba^{i_1+1}b\bigr)+\cdots  \cr
&\ \ \ \ \ \ +\bigl(f\,|\,a^{i_r+1}ba^{i_1}b\dots ba^{i_{r-1}}b\bigr)-
\bigl(f\,|\,a^{i_r}ba^{i_1}b\cdots ba^{i_{r-1}+1}b\bigr)\cr
&=-Kr\,\Bigl([a,b]^r\,|\,a^{i_1}b\cdots a^{i_r}b\Bigr)+(-1)^{r-1}
Kr\,\Bigl([a,b]^r\,|\,a^{i_r}b\cdots a^{i_1}b\Bigr).}\eqno(14)$$
Since $f$ is push-invariant, we have $(f|ub)=(f|bu)$ for every word $u$,
so we can modify the negative terms in (14):
$$\eqalign{&\bigl(f\,|\,a^{i_1+1}ba^{i_2}b\dots ba^{i_r}b\bigr)-
\bigl(f\,|\,ba^{i_1}ba^{i_2}b\cdots ba^{i_r+1}\bigr)  \cr
&\ \ \ +\bigl(f\,|\,a^{i_2+1}ba^{i_3}b\dots ba^{i_1}b\bigr)-
\bigl(f\,|\,ba^{i_2}ba^{i_3}b\cdots ba^{i_1+1}\bigr)+\cdots  \cr
&\ \ \ \ \ \ +\bigl(f\,|\,a^{i_r+1}ba^{i_1}b\dots ba^{i_{r-1}}b\bigr)-
\bigl(f\,|\,ba^{i_r}ba^{i_1}b\cdots ba^{i_{r-1}+1}\bigr)\cr
&=-Kr\,\Bigl([a,b]^r\,|\,a^{i_1}b\cdots a^{i_r}b\Bigr)+(-1)^{r-1}
Kr\,\Bigl([a,b]^r\,|\,a^{i_r}b\cdots a^{i_1}b\Bigr).}\eqno(15)$$
Now all words in the positive terms start in $a$ and end in $b$, and
all words int he negative terms start in $b$ and end in $a$, so
we can remove these letters and write
$$\eqalign{&\bigl(f^a_b\,|\,a^{i_1}ba^{i_2}b\dots ba^{i_r}\bigr)-
\bigl(f^b_a\,|\,a^{i_1}ba^{i_2}b\cdots ba^{i_r}\bigr)  \cr
&\ \ \ +\bigl(f^a_b\,|\,a^{i_2}ba^{i_3}b\dots ba^{i_1}\bigr)-
\bigl(f^b_a\,|\,a^{i_2}ba^{i_3}b\cdots ba^{i_1}\bigr)+\cdots  \cr
&\ \ \ \ \ \ +\bigl(f^a_b\,|\,a^{i_r}ba^{i_1}b\dots ba^{i_{r-1}}\bigr)-
\bigl(f^b_a\,|\,a^{i_r}ba^{i_1}b\cdots ba^{i_{r-1}}\bigr)\cr
&=-Kr\,\Bigl([a,b]^r\,|\,a^{i_1}b\cdots a^{i_r}b\Bigr)+
(-1)^{r-1}Kr\,\Bigl([a,b]^r\,|\,a^{i_r}b\cdots a^{i_1}b\Bigr).}\eqno(16)$$
The left-hand side of this equal to
$$\bigl(pushsym(f^a_b-f^b_a)\,|\,a^{i_1}b\cdots ba^{i_r}\bigr),$$
so to the left-hand side of the divergence condition (5) for the weight
$n-2$ word $u=a^{i_1}b\cdots ba^{i_r}$.  If $n\ne 2r+1$, then the right-hand 
sides of (5) and (16) are both equal to $0$.  To prove that (5) and (16) are
identical, it remains only to check that the right-hand 
sides are equal when $n=2r+1$, which, cancelling the factor $Kr$ from
both sides, reduces to the following lemma. 
\vskip .3cm
\noindent {\bf Lemma.} {\it For each word $u$ of depth $r-1$ and weight $2r-1$, 
we have
$${{1}\over{|C_b(ub)|}}\sum_{v\in C(ub)} \bigl([a,b]^r\,|\,v\bigr)=
\bigl([a,b]^r\,|\,ub\bigr)-(-1)^{r-1}
\bigl([a,b]^r\,|\,u'b\bigr),\eqno(17)$$
where $u'$ denotes the word $u$ written backwards.}
\vskip .3cm
\noindent Proof.  Observe that if $([a,b]^r|ub)\ne 0$, then $ub$ must satisfy 
the {\it parity property} that, writing $ub=u_1\cdots u_{2r}$ where each $u_i$ 
is letter $a$ or $b$, the pair $u_{2i-1}u_{2i}$ must be either $ab$ or 
$ba$ for $0\le i\le r$.  The coefficient of the word $ub$ in $[a,b]^r$ is 
equal to $(-1)^j$ where $j$ is the number of pairs $u_{2i-1}u_{2i}$ in $ub$ 
that are equal to $ba$.  In other words, if a word $w$ appears with non-zero
coefficient in $[a,b]^r$, then letting $U=ba$ and $V=ab$, we must be able to 
write $w$ as a word in $U,V$, and the coefficient of $w$ in $[a,b]^r$ is
$(-1)^m$ where $m$ denotes the number of times the letter $U$ occurs.

If $w=ub=V^r=(ab)^r$, then $u'b=ub$. The coefficient of $V^r$ in $[a,b]^r$
is equal to 1, so the right-hand side of (17) is
equal to 2 if $r$ is even and $0$ if $r$ is odd. For the left-hand side,
$C(ub)=\{V^r,U^r\}$ and $C_b(ub)=\{V^r\}$, so $|C_b(ub)|=1$.  The coefficient of
$U^r$ in $[a,b]^r$ is equal to $(-1)^r$, so the left-hand side is again 
equal to 2 if $r$ is even and $0$ if $r$ is odd.  This proves (17) in the
case $ub=V^r$.

Suppose now that $ub\ne V^r$ but that it satisfies the parity property. 
Write $ub=U^{a_1}V^{b_1}\cdots U^{a_s}V^{b_s}$ in which all the $a_i,b_i\ge 1$ 
except for $a_1$, which may be $0$.  Then $u'b$ is equal to 
$aU^{b_s-1}V^{a_s}\cdots U^{b_1}V^{a_1}b$. If $b_s>1$, then the
pair $u_{2(b_s-1)+1}u_{2(b_s-1)+2}$ is $aa$, so $([a,b]^r|u'b)=0$. 
If $b_s=1$, then the word $u'b$ begins with $aa$ and thus does not have
the parity property, so again $([a,b]^r|u'b)=0$.
This shows that if $([a,b]^r|ub)\ne 0$ then $([a,b]^r|u'b)=0$ and
vice versa. 

This leaves us with three possibilities for $ub\ne V^r$.
\vskip .1cm
\noindent {\it Case 1: $([a,b]^r|ub)\ne 0$.} Then $ub$ has the parity
property, so we write $ub=U^{a_1}V^{b_1}\cdots U^{a_s}V^{a_s}$ as above.
The right-hand side of (17) is then equal to $(-1)^j$ where $j=a_1+\cdots+a_s$.
For the left-hand side, we
note that the only words in the cyclic permutation class $C(ub)$ that have the
parity property are the cyclic shifts of $ub$ by an even number of letters,
otherwise a pair $aa$ or $bb$ necessarily occurs as above. These are the
same as the cyclic permutations of the word $ub$ written in the letters $U,V$.  
All these cyclic permutations obviously have the same number of occurrences 
$j$ of the letter $U$. Thus, the words in $C(ub)$ for which $[a,b]^r$ has a 
non-zero coefficient are the cyclic permutations of the word $ub$ in the
letters $U,V$, and the coefficient is always equal to $(-1)^j$.  
These words are exactly half of the all the words in $C(ub)$, 
so the sum in the left-hand side is equal to 
$(-1)^j|C(ub)|/2$.  But $|C_b(ub)|=|C(ub)|/2$, so the left-hand side is
equal to $(-1)^j$, which proves (17) for words $ub$ having the parity property.
\vskip .1cm
\noindent {\it Case 2: $([a,b]^r|u'b)\ne 0$.} In this case it is $u'b$ that
has the parity property, and the right-hand side of (17) is equal to
$(-1)^{r+j'}$ where $j'$ is the number of occurrences of $U$ in the word
$u'b=U^{a_1}V^{b_1}\cdots U^{a_s}V^{b_s}$, i.e.~$j'=a_1+\cdots+a_s$. We have
$ub=aU^{b_s-1}V^{a_s}\cdots U^{b_1}V^{b_1}b$. The word
$w=U^{b_s-1}V^{a_s}\cdots  U^{b_1}V^{b_1}U$ then occurs in $C(ub)$, and
the number of occurrences of the letter $U$ in $ub$ is equal to
$j=b_1+\cdots+b_{s-1}+b_s$. Since $a_1+b_1+\cdots+a_s+b_s=r$, we
have $j+j'=r$ so $j'=r-j$ and the right-hand side of (17) is equal to
$(-1)^j$. The number of words in $C(ub)$ that have non-zero coefficient in
$[a,b]^r$ is $|C(ub)|/2=|C_b(ub)|$ as above, these words being exactly the
cyclic permutations of $w$ written in $U,V$, and the coefficient is always
equal to $(-1)^j$.  So the left-hand side of (17) is equal to $(-1)^j$, 
which proves (17) in the case where $u'b$ has the parity property.
\vskip .2cm
\noindent {\it Case 3: $([a,b]^r|ub)=([a,b]^r|u'b)=0$.} The right-hand
side of (17) is zero.  For the left-hand side, consider the words in
$C(ub)$.  If there are no words in $C(ub)$ whose coefficient in
$[a,b]^r$ is non-zero, then the left-hand side of (17) is also zero and 
(17) holds.  Suppose instead that
there is a word $w\in C(ub)$ whose coefficient in $[a,b]^r$ is non-zero. Then
as we saw above, $w$ is a cyclic shift of $ub$ by an odd number of letters,
and since all cyclic shifts of $w$ by an even number of letters then have
the same coefficient in $[a,b]^r$ as $w$, we may assume that $w$ is the
cyclic shift of $ub$ by one letter, i.e.~taking the final $b$ and putting 
it at the beginning.  Since $w$ has non-zero coefficient in $[a,b]^r$,
we can write $w=U^{a_1}V^{b_1}\cdots U^{a_s} V^{b_s}$,  where $a_1>0$ since
$w$ now starts with $b$, but $b_s$ may be equal to $0$ since $w$ may end 
with $a$.  Then $ub=aU^{a_1-1}V^{b_1}\cdots U^{a_s}V^{b_s}b$, so
we can write $u'b=U^{b_s}V^{a_s}\cdots U^{b_1}V^{a_1-1}ab=U^{b_s}V^{a_s}
\cdots U^{v_1}V^{a_1}$. But then $u'b$ satisfies the parity property, so
its coefficient in $[a,b]^r$ is non-zero, contradicting our assumption.
Thus under the assumption, all words in $C(ub)$ have coefficient zero in
$[a,b]^r$, which completes the proof of the Lemma.\hfill{$\diamondsuit$}
\vskip 1.3cm
\noindent {\bf References}
\vskip .4cm
\noindent [AKKN1] A.~Alekseev, N.~Kawazumi, Y.~Kuno, F.~Naef, 
Higher genus Kashiwara-Vergne problems and the Goldman-Turaev Lie
bialgebra, arXiv:1611.05581, November 2016.
\vskip .2cm
\noindent [AKKN2] A.~Alekseev, N.~Kawazumi, Y.~Kuno, F.~Naef, The
Goldman-Turaev Lie bialgebra and the Kashiwara-Vergne problem in higher
genera, arXiv:1804.09566v2, May 2018.
\vskip .2cm
\noindent [AT] A.~Alekseev, C.~Torossian, The Kashiwara-Vergne conjecture
and Drinfel'd's associators, {\it Annals of Math.} {\bf 175} no. 2 (2012), 
415-463.
\vskip .2cm
\noindent [CEE] D.~Calaque, P.~Etingof, B.~Enriquez, Universal KZB equations: 
the elliptic case, {\it Algebra, Arithmetic and Geometry: in honor of 
Yu.~I.~Manin}, Vol. I, 165-266, {\it Progr. Math.} {\bf 269}, Birkh\"auser 
Boston, 2009.
\vskip .2cm
\noindent [D] V.G.~Drinfel'd, On quasitriangular quasi-Hopf algebras and
a group that is closely connected with ${\rm Gal}(\overline{\Bbb Q}/
{\Bbb Q})$, {\it Algebra i Analiz} {\bf 2} (1990), 149-181; translation
in {\it Leningrad Math. J.} {\bf 2} (1991), 829-860.
\vskip .2cm
\noindent [E] J.~\'Ecalle, The flexion structure and dimorphy: flexion units, 
singulators, generators, and the enumeration of multizeta irreducibles, 
in {\it Asymptotics in Dynamics, Geometry and PDEs; Generalized Borel 
Summation vol. II}, Costin O., Fauvet F., Menous F., Sauzin D., eds., 
CRM Series, vol. 12.2. Edizioni della Normale, 2011.
\vskip .2cm
\noindent [En] B.~Enriquez, Elliptic associators, {\it Selecta Math. New
Series} {\bf 20} (2014), 491-584.
\vskip .2cm
\noindent [F] H.~Furusho, Double shuffle relation for associators,
{\it Annals of Math.} {\bf 174} no. 1 (2011), 341-360.
\vskip .2cm
\noindent [G] A.~Grothendieck, Esquisse d'un Programme, in {\it Geometric
Galois Actions 1}, eds. L.~Schneps, P.~Lochak, London Math. Soc. Lecture
Notes {\bf 242} (1997), 5-48.
\vfill\eject
\noindent [H] R.~Hain, Infinitesimal presentations of the Torelli groups,
{\it J. AMS} {\bf 10} no. 3 (1997), 597-651.
\vskip .2cm
\noindent [HLS] A.~Hatcher, P.~Lochak, L.~Schneps, On the Teichm\"uller
tower of mapping class groups, {\it J. reine angew. Math.} {\bf 521} (2000),
1-24.
\vskip .2cm
\noindent [I1] Y.~Ihara, Braids, Galois groups, and some arithmetic functions,
in {\it Proceedings of the ICM, Kyoto, Japan 1990}, Math. Soc. Japan, 1991.
\vskip .2cm
\noindent [I2] Y.~Ihara, On the stable derivation algebra associated
with some braid groups, {\it Israel J. Math.} {\bf 80} no. 1-2 (1992),
135-153.
\vskip .2cm
\noindent [NS] H.~Nakamura, L.~Schneps, On a subgroup of the 
Grothendieck-Teichm\"uller group acting on the tower of profinite
Teichm\"uller modular groups, {\it Invent. Math.} {\bf 141} no. 3 (2000), 
503-560.
\vskip .2cm
\noindent [RS] E.~Raphael, L.~Schneps, On linearised and elliptic versions
of the Kashiwara-Vergne Lie algebra, arXiv:1706.08299, June 2017.
\vskip .2cm
\noindent [S1] L.~Schneps, Double shuffle and Kashiwara-Vergne Lie algebras,
{\it J. Algebra} {\bf 367} (2012), 54-74.
\vskip .2cm
\noindent [S2] L.~Schneps, ARI, GARI, Zig and Zag: An
introduction to \'Ecalle's theory of multiple zeta values, arXiv:1507.01534,
July 2015.
\vskip .2cm
\noindent [S3] L.~Schneps, Elliptic double shuffle,
Grothendieck-Teichm\"uller and mould theory, arXiv:1506.09050, June 2015.
\vskip .2cm
\noindent [T] H.~Tsunogai, The stable derivation algebras for higher
genera, {\it Israel J. Math.} {\bf 136} (2003), 221-250.
\vskip 1cm
{\tiny \'Elise Raphael 

\vskip -.1cm
Universit\'e de Gen\`eve 

\vskip -.1cm
Elise.Raphael@unige.ch}
\vskip .2cm
{\tiny Leila Schneps 

\vskip -.1cm
Institut de Math\'ematiques de Jussieu 

\vskip -.1cm
Leila.Schneps@imj-prg.fr}

\bye